\documentclass{gtart}

\def\ifplaintex{\expandafter\ifx\csname documentclass\endcsname\relax}

\ifplaintex 
\else
\fi

\expandafter\ifx\csname beginpicture\endcsname\relax
\expandafter\ifx\csname documentclass\endcsname\relax
\input pictex \else
\input prepictex \input pictex \input postpictex \fi\fi

\def\gt{{\mathsurround=0pt\it $\cal G\mskip-2mu$eometry \&\ 
$\cal T\!\!$opology}}        

\def\gtp{{\mathsurround=0pt\it $\cal G\mskip-2mu$eometry \&\ 
$\cal T\!\!$opology $\cal P\!$ublications}}  

\def\lognumber#1{\def\thelognumber{#1}}
\def\volumenumber#1{\def\thevolumenumber{#1}}
\def\papernumber#1{\def\thepapernumber{#1}}
\def\volumeyear#1{\def\thevolumeyear{#1}}

\def\pagenumbers#1#2{\def\startpage{#1}\def\finishpage{#2}}
\def\published#1{\def\publishdate{#1}}
\def\proposed#1{\def\theproposer{#1}}
\def\seconded#1{\def\theseconders{#1}}
\def\received#1{\def\receiveddate{#1}}
\def\revised#1{\def\reviseddate{#1}}
\def\accepted#1{\def\accepteddate{#1}}
\def\asciititle#1{\def\theasciititle{#1}}
\def\covertitle#1{\def\thecovertitle{#1}}

\long\def\asciiabstract#1{\long\def\theasciiabstract{#1}}
\def\asciikeywords#1{\def\theasciikeywords{#1}}

\def\shorttitle#1{\def\theshorttitle{#1}}

\let\\\par\let\thelognumber\relax
\let\thevolumenumber\relax\let\thepapernumber\relax
\let\thevolumeyear\relax\let\thesamplenumber\relax\let\startpage\relax
\let\finishpage\relax\let\publishdate\relax\let\receiveddate\relax
\let\reviseddate\relax\let\accepteddate\relax\let\theasciititle\relax
\let\thecovertitle\relax\let\theasciiauthors\relax
\let\theasciiabstract\relax\let\theasciikeywords\relax
\let\theasciiemail\relax\let\theshortauthors\relax\let\theshorttitle\relax

\long\def\maketitlep{   

\count0=\startpage

\gt\hfill      
\beginpicture
\setcoordinatesystem units <0.33truein, 0.33truein> point at 2.2 0.9
\setplotsymbol ({$\cal G$})
\plotsymbolspacing=9truept
\circulararc 315 degrees from 0 1 center at 0 0
\setplotsymbol ({$\cal T$})
\circulararc 315 degrees from 1 -1 center at 1 0
\endpicture
\break
{\small\ifx\thesamplenumber\relax 
Volume \else Sample
\fi\thevolumenumber\ (\thevolumeyear)
\startpage--\finishpage\nl
Published: \publishdate}
\vglue 0.5truein plus 0.4fil minus 0.1truein

{\parskip=0pt\leftskip 0pt plus 1fil\def\\{\par\smallskip}{\ifplaintex\large
\else\Large\fi\bf\thetitle}\par\medskip}   

\vglue 0pt plus 0.1fil 

{\parskip=0pt\leftskip 0pt plus 1fil\def\\{\par}{\sc\theauthors}
\par\medskip}

\vglue 0pt plus 0.1fil 

{\small\parskip=0pt\let\newline\\
{\leftskip 0pt plus 1fil\def\\{\par}{\sl\theaddress}\par}
\expandafter\ifx\theemail\relax    
\relax\else\vglue 5pt plus 0.02fil minus 2pt\def\\{\stdspace{\rm 
and}\stdspace} 
\cl{Email:\stdspace\tt\theemail}\fi
\ifx\theurl\relax                  
\relax\else\vglue 5pt plus 0.02fil minus 2pt\def\\{\stdspace{\rm 
and}\stdspace}
\cl{URL:\stdspace\tt\theurl}\fi\par}

\vglue 7pt plus 0.3fil minus 3pt

{\bf Abstract}
\vglue 5pt plus 0.1fil minus 2pt

\theabstract

\vglue 7pt plus 0.3fil minus 3pt

{\bf AMS Classification numbers}\quad Primary:\quad \theprimaryclass

Secondary:\quad \thesecondaryclass

\vglue 5pt plus 0.3fil minus 2pt

{\bf Keywords}\quad \thekeywords

\vglue 10pt plus 0.5fil minus 5pt

{\small  Proposed: \theproposer\hfill Received: \receiveddate\nl
Seconded: \theseconders\hfill 
\ifx\reviseddate\relax                         
Accepted: \accepteddate                        
\else
Revised: \reviseddate                          
\fi}
\eject
}       

\let\maketitlepage\maketitlep
\let\maketitle\maketitlepage

\font\phead=cmsl9 scaled 950
\font=cmsl9 scaled 1050
\font\pnum=cmbx10 scaled 913
\font\lnum=cmbx10 
\font\pfoot=cmsl9 scaled 950
\font=cmsl9 scaled 1050
\ifplaintex
\headline{\vbox to 0pt{\vskip -4.5mm\line{\small\phead\ifnum
\count0=\startpage ISSN 1364-0380 (on line)
1465-3060 (printed) \hfill {\pnum\folio}\else\ifodd\count0\def\\{ }%
\ifx\theshorttitle\relax\thetitle\else\theshorttitle\fi\hfill{\pnum\folio}
\else\def\\{ and }{\pnum\folio}\hfill\ifx\theshortauthors\relax\theauthors
\else\theshortauthors\fi\fi\fi}\vss}}
\footline{\vbox to 0pt{\vglue 0mm\line{\small\pfoot\ifnum\count0=\startpage
\copyright\ \gtp\hfill\else
\gt, Volume \thevolumenumber\ (\thevolumeyear)\hfill\fi}\vss
}}
\else
\makeatletter
\def\@oddhead{{\small\lhead\ifnum\count0=\startpage ISSN 1364-0380 (on line)
1465-3060 (printed) \hfill {\lnum\number\count0}\else\ifodd\count0
\def\\{ }\ifx\theshorttitle\relax \thetitle \else\theshorttitle\fi\hfill
{\lnum\number\count0}\else\def\\{ and }{\lnum\number\count0}
\hfill\ifx\theshortauthors\relax 
\theauthors\else\theshortauthors\fi\fi\fi}}\def\@evenhead{\@oddhead}
\def\@oddfoot{\small\lfoot\ifnum\count0=\startpage\copyright\ \gtp\hfill\else
\gt, Volume \thevolumenumber\ (\thevolumeyear)\hfill\fi}
\def\@evenfoot{\@oddfoot}
\makeatother
\fi

\newwrite\gtoutfile
\long\gdef\makeheadfile{  
{\def\\{, }\def\s{ }
\immediate\openout\gtoutfile head.xxx
\immediate\write\gtoutfile{To: math@arxiv.org}
\immediate\write\gtoutfile{Subject: put or rep NNNNN:pppp}
\immediate\write\gtoutfile{--text follows this line--}
\immediate\write\gtoutfile{Proxy-for: \ifx\theasciiauthors\relax
\theauthors\else\theasciiauthors\fi\s<\ifx\theasciiemail\relax\theemail\else\theasciiemail\fi>}
\immediate\write\gtoutfile{\noexpand\\}
\immediate\write\gtoutfile{Authors: \ifx\theasciiauthors\relax
\theauthors\else\theasciiauthors\fi}
{\def\\{ }\immediate\write\gtoutfile{Title: \ifx\theasciititle\relax
\thetitle\else\theasciititle\fi}}
\immediate\write\gtoutfile{Subj-class: GT or SG or MG etc}
\immediate\write\gtoutfile{MSC-class: \theprimaryclass\ifx\thesecondaryclass\relax\else, \thesecondaryclass\fi}
\immediate\write\gtoutfile{Journal-ref: Geom. Topol. \thevolumenumber
(\thevolumeyear) \startpage-\finishpage}
\immediate\write\gtoutfile{Comments: Published by Geometry and Topology at}
\immediate\write\gtoutfile{\s\s http://www.maths.warwick.ac.uk/gt/GTVol\thevolumenumber/paper\thepapernumber.abs.html}
\immediate\write\gtoutfile{\noexpand\\}
\immediate\write\gtoutfile{}
\ifx\theasciiabstract\relax
\immediate\write\gtoutfile{\theabstract}\else
\immediate\write\gtoutfile{\theasciiabstract}\fi
\immediate\write\gtoutfile{}
\immediate\write\gtoutfile{\noexpand\\}
\immediate\write\gtoutfile{}
\immediate\closeout\gtoutfile}}  

\def\maketitlepage{\maketitlep\makeheadfile}
\let\maketitle\maketitlepage

\def\ifplaintex{\expandafter\ifx\csname documentclass\endcsname\relax}

\ifplaintex 
\else
\fi

\expandafter\ifx\csname beginpicture\endcsname\relax
\expandafter\ifx\csname documentclass\endcsname\relax
\input pictex \else
\input prepictex \input pictex \input postpictex \fi\fi

\def\gt{{\mathsurround=0pt\it $\cal G\mskip-2mu$eometry \&\ 
$\cal T\!\!$opology}}        

\def\gtp{{\mathsurround=0pt\it $\cal G\mskip-2mu$eometry \&\ 
$\cal T\!\!$opology $\cal P\!$ublications}}  

\def\lognumber#1{\def\thelognumber{#1}}
\def\volumenumber#1{\def\thevolumenumber{#1}}
\def\papernumber#1{\def\thepapernumber{#1}}
\def\volumeyear#1{\def\thevolumeyear{#1}}

\def\pagenumbers#1#2{\def\startpage{#1}\def\finishpage{#2}}
\def\published#1{\def\publishdate{#1}}
\def\proposed#1{\def\theproposer{#1}}
\def\seconded#1{\def\theseconders{#1}}
\def\received#1{\def\receiveddate{#1}}
\def\revised#1{\def\reviseddate{#1}}
\def\accepted#1{\def\accepteddate{#1}}
\def\asciititle#1{\def\theasciititle{#1}}
\def\covertitle#1{\def\thecovertitle{#1}}

\long\def\asciiabstract#1{\long\def\theasciiabstract{#1}}
\def\asciikeywords#1{\def\theasciikeywords{#1}}

\def\shorttitle#1{\def\theshorttitle{#1}}

\let\\\par\let\thelognumber\relax
\let\thevolumenumber\relax\let\thepapernumber\relax
\let\thevolumeyear\relax\let\thesamplenumber\relax\let\startpage\relax
\let\finishpage\relax\let\publishdate\relax\let\receiveddate\relax
\let\reviseddate\relax\let\accepteddate\relax\let\theasciititle\relax
\let\thecovertitle\relax\let\theasciiauthors\relax
\let\theasciiabstract\relax\let\theasciikeywords\relax
\let\theasciiemail\relax\let\theshortauthors\relax\let\theshorttitle\relax

\long\def\maketitlep{   

\count0=\startpage

\gt\hfill      
\beginpicture
\setcoordinatesystem units <0.33truein, 0.33truein> point at 2.2 0.9
\setplotsymbol ({$\cal G$})
\plotsymbolspacing=9truept
\circulararc 315 degrees from 0 1 center at 0 0
\setplotsymbol ({$\cal T$})
\circulararc 315 degrees from 1 -1 center at 1 0
\endpicture
\break
{\small\ifx\thesamplenumber\relax 
Volume \else Sample
\fi\thevolumenumber\ (\thevolumeyear)
\startpage--\finishpage\nl
Published: \publishdate}
\vglue 0.5truein plus 0.4fil minus 0.1truein

{\parskip=0pt\leftskip 0pt plus 1fil\def\\{\par\smallskip}{\ifplaintex\large
\else\Large\fi\bf\thetitle}\par\medskip}   

\vglue 0pt plus 0.1fil 

{\parskip=0pt\leftskip 0pt plus 1fil\def\\{\par}{\sc\theauthors}
\par\medskip}

\vglue 0pt plus 0.1fil 

{\small\parskip=0pt\let\newline\\
{\leftskip 0pt plus 1fil\def\\{\par}{\sl\theaddress}\par}
\expandafter\ifx\theemail\relax    
\relax\else\vglue 5pt plus 0.02fil minus 2pt\def\\{\stdspace{\rm 
and}\stdspace} 
\cl{Email:\stdspace\tt\theemail}\fi
\ifx\theurl\relax                  
\relax\else\vglue 5pt plus 0.02fil minus 2pt\def\\{\stdspace{\rm 
and}\stdspace}
\cl{URL:\stdspace\tt\theurl}\fi\par}

\vglue 7pt plus 0.3fil minus 3pt

{\bf Abstract}
\vglue 5pt plus 0.1fil minus 2pt

\theabstract

\vglue 7pt plus 0.3fil minus 3pt

{\bf AMS Classification numbers}\quad Primary:\quad \theprimaryclass

Secondary:\quad \thesecondaryclass

\vglue 5pt plus 0.3fil minus 2pt

{\bf Keywords}\quad \thekeywords

\vglue 10pt plus 0.5fil minus 5pt

{\small  Proposed: \theproposer\hfill Received: \receiveddate\nl
Seconded: \theseconders\hfill 
\ifx\reviseddate\relax                         
Accepted: \accepteddate                        
\else
Revised: \reviseddate                          
\fi}
\eject
}       

\let\maketitlepage\maketitlep
\let\maketitle\maketitlepage

\font\phead=cmsl9 scaled 950
\font=cmsl9 scaled 1050
\font\pnum=cmbx10 scaled 913
\font\lnum=cmbx10 
\font\pfoot=cmsl9 scaled 950
\font=cmsl9 scaled 1050
\ifplaintex
\headline{\vbox to 0pt{\vskip -4.5mm\line{\small\phead\ifnum
\count0=\startpage ISSN 1364-0380 (on line)
1465-3060 (printed) \hfill {\pnum\folio}\else\ifodd\count0\def\\{ }%
\ifx\theshorttitle\relax\thetitle\else\theshorttitle\fi\hfill{\pnum\folio}
\else\def\\{ and }{\pnum\folio}\hfill\ifx\theshortauthors\relax\theauthors
\else\theshortauthors\fi\fi\fi}\vss}}
\footline{\vbox to 0pt{\vglue 0mm\line{\small\pfoot\ifnum\count0=\startpage
\copyright\ \gtp\hfill\else
\gt, Volume \thevolumenumber\ (\thevolumeyear)\hfill\fi}\vss
}}
\else
\makeatletter
\def\@oddhead{{\small\lhead\ifnum\count0=\startpage ISSN 1364-0380 (on line)
1465-3060 (printed) \hfill {\lnum\number\count0}\else\ifodd\count0
\def\\{ }\ifx\theshorttitle\relax \thetitle \else\theshorttitle\fi\hfill
{\lnum\number\count0}\else\def\\{ and }{\lnum\number\count0}
\hfill\ifx\theshortauthors\relax 
\theauthors\else\theshortauthors\fi\fi\fi}}\def\@evenhead{\@oddhead}
\def\@oddfoot{\small\lfoot\ifnum\count0=\startpage\copyright\ \gtp\hfill\else
\gt, Volume \thevolumenumber\ (\thevolumeyear)\hfill\fi}
\def\@evenfoot{\@oddfoot}
\makeatother
\fi

\newwrite\gtoutfile
\long\gdef\makeheadfile{  
{\def\\{, }\def\s{ }
\immediate\openout\gtoutfile head.xxx
\immediate\write\gtoutfile{To: math@arxiv.org}
\immediate\write\gtoutfile{Subject: put or rep NNNNN:pppp}
\immediate\write\gtoutfile{--text follows this line--}
\immediate\write\gtoutfile{Proxy-for: \ifx\theasciiauthors\relax
\theauthors\else\theasciiauthors\fi\s<\ifx\theasciiemail\relax\theemail\else\theasciiemail\fi>}
\immediate\write\gtoutfile{\noexpand\\}
\immediate\write\gtoutfile{Authors: \ifx\theasciiauthors\relax
\theauthors\else\theasciiauthors\fi}
{\def\\{ }\immediate\write\gtoutfile{Title: \ifx\theasciititle\relax
\thetitle\else\theasciititle\fi}}
\immediate\write\gtoutfile{Subj-class: GT or SG or MG etc}
\immediate\write\gtoutfile{MSC-class: \theprimaryclass\ifx\thesecondaryclass\relax\else, \thesecondaryclass\fi}
\immediate\write\gtoutfile{Journal-ref: Geom. Topol. \thevolumenumber
(\thevolumeyear) \startpage-\finishpage}
\immediate\write\gtoutfile{Comments: Published by Geometry and Topology at}
\immediate\write\gtoutfile{\s\s http://www.maths.warwick.ac.uk/gt/GTVol\thevolumenumber/paper\thepapernumber.abs.html}
\immediate\write\gtoutfile{\noexpand\\}
\immediate\write\gtoutfile{}
\ifx\theasciiabstract\relax
\immediate\write\gtoutfile{\theabstract}\else
\immediate\write\gtoutfile{\theasciiabstract}\fi
\immediate\write\gtoutfile{}
\immediate\write\gtoutfile{\noexpand\\}
\immediate\write\gtoutfile{}
\immediate\closeout\gtoutfile}}  

\def\maketitlepage{\maketitlep\makeheadfile}
\let\maketitle\maketitlepage

\lognumber{166}
\volumenumber{6}\papernumber{1}\volumeyear{2002}
\pagenumbers{1}{26}
\accepted{12 January 2002}
\proposed{David Gabai}
\seconded{Jean-Pierre Otal, Robion Kirby}
\received{20 February 2001}
\revised{26 October 2001}
\published{16 January 2002}

\usepackage{amssymb,graphicx}

\def\R{\mathbb{R}}
\def\Z{\mathbb{Z}}
\def\C{\mathbb{C}}
\def\Q{\mathbb{Q}}
\def\H{\mathbb{H}}
\def\B{\mathbb{B}}

\newtheorem{theorem}{Theorem}[section]
\newtheorem{lemma}[theorem]{Lemma}

\theoremstyle{definition}
\newtheorem{remark}[theorem]{Remark}
\newtheorem{definition}[theorem]{Definition}

\newcommand{\pim}{\ensuremath{\pi_1 M}}
\newcommand{\MARGEST}{\ensuremath{\frac{.49}{16(2(4\pi/.49)^3+1)}}}
\newcommand{\Qbar}{\ensuremath{\overline{\mathbb{Q}}}}
\newenvironment{Relax}{\relax}{\relax}
\begin{document}

\title{Algorithmic detection and description of\\\vglue -5pt\\hyperbolic structures
on closed 3--manifolds\\\vglue -1pt\\with solvable word problem}
\covertitle{Algorithmic detection and description of\\hyperbolic structures
on closed 3--manifolds\\with solvable word problem}
\asciititle{Algorithmic detection and description of hyperbolic structures
on closed 3-manifolds with solvable word problem}
\shorttitle{Algorithmic detection of hyperbolic structures}

\author{Jason Manning}                  

\address{Department of Mathematics, University of California at 
Santa Barbara\\Santa Barbara, CA 93106, USA}
\email{manning@math.ucsb.edu}                     

\begin{abstract} 
We outline a rigorous algorithm, first suggested
by Casson, for determining whether a closed orientable $3$-manifold $M$ 
is hyperbolic, and to compute the hyperbolic structure, if one exists.
The algorithm requires that a procedure has been given to solve the
word problem in $\pim$. 
\end{abstract}
\begin{Relax}\end{Relax}

\asciiabstract{ 
We outline a rigorous algorithm, first suggested
by Casson, for determining whether a closed orientable 3-manifold M 
is hyperbolic, and to compute the hyperbolic structure, if one exists.
The algorithm requires that a procedure has been given to solve the
word problem in \pi_1(M).}

\asciikeywords{3-manifold, Kleinian group, word problem, recognition problem,
geometric structure}                    

\maketitlepage

\section{Introduction}\label{sec:intro}

Unless otherwise noted, all manifolds are assumed to be orientable. 
We begin with a couple of definitions:
\begin{definition}
We say that a $3$--manifold $M$ has a \emph{geometric structure} (or 
that $M$ is \emph{geometric}) if there is a
homeomorphism 
$f\co \mathrm{Interior}(M)\to N$ where $N$ is a complete homogeneous
Riemannian manifold which 
is locally isometric to one of the eight model geometries,
which are discussed in \cite{scott:geom}.  If $N$ is locally isometric to
hyperbolic $3$--space, we say $M$ has a \emph{hyperbolic structure}.  
The terms \emph{spherical structure}, \emph{Euclidean structure}, and so
on are similarly defined.

If there is an algorithm to construct $N$ (for instance as a finite number of
charts or as a convex polyhedron in the model geometry with face 
identifications) and an algorithm to
construct the homeomorphism $f\co M \to N$, we say that the geometric
structure on $M$ is \emph{algorithmically constructible}.
\end{definition}

\begin{remark}\label{remark:mapfollows}
If two triangulated compact $3$--manifolds $M$ and $N$ are homeomorphic, 
then it is shown in \cite{moise:hauptvermutung} (see also \cite{bing:hauptvermutung}) that the given 
triangulations admit isomorphic
subdivisions. 
For any $n$, there are finitely many subdivisions of the
triangulations of $N$ and $M$ with $n$ tetrahedra, so we can find them all and
test each pair for isomorphism.  Eventually we must find subdivisions which
are isomorphic, and the isomorphism gives a homeomorphism 
between the two manifolds.
Thus to prove that the geometric structure on a triangulated
$3$--manifold $M$ is algorithmically constructible it is necessary only to
produce a triangulation of $N$.  
\end{remark}

\begin{definition}
Let $p$ be some statement.  An algorithm 
\emph{decides if
$p$} if whenever $p$ is true, the algorithm reports that $p$ is true, and
whenever $p$ is false, the algorithm remains silent.  
An algorithm \emph{decides whether or not}
$p$ if whenever $p$ is true the algorithm reports that $p$ is true, and 
whenever $p$ is false the
algorithm reports that $p$ is false.
\end{definition}

It is an open question whether there is an algorithm to decide whether 
or not a closed triangulated $3$--manifold has a hyperbolic structure.  
In this paper we prove:
\medskip

\noindent
{\bf Theorem \ref{th:main}}\qua {\sl
There exists an algorithm which will, given a triangulated orientable closed 
$3$-manifold $M$ and a
solution to the word problem in $\pim$, decide whether or not $M$ has a
hyperbolic structure.
}

\medskip

The algorithm described gives us detailed information about the
hyperbolic structure, and we are able to show further that:

\medskip

\noindent
{\bf Theorem \ref{th:structure}}\qua {\sl
If $M$ is a triangulated closed orientable $3$--manifold which has a
hyperbolic structure, then the hyperbolic structure is algorithmically
constructible.
}

\medskip

\subsection{Context}

Several classes of manifolds 
are
known to have fundamental groups with solvable word problems;  Waldhausen 
solved the word problem for
Haken and virtually Haken manifolds \cite{waldhausen:wo}, Niblo and Reeves for
manifolds (of any dimension) with a 
nonpositively curved cubing \cite{nib:cube}.  
Skinner has
solved the word problem in a class of non-Haken manifolds containing
immersed surfaces of a particular type \cite{skinner:word}.  Manifolds with
word hyperbolic fundamental group also have solvable word problem 
\cite{ep:word}.  Lackenby \cite{lackenby:whds} has shown that there are many
$3$-manifolds with word hyperbolic fundamental group, not all of which are known to be
hyperbolic.  If Thurston's geometrization conjecture \cite{thurston:bams} is 
true, then these
manifolds are in fact hyperbolic.

If $M$ is any $3$--manifold satisfying the geometrization conjecture, a solution to 
the word problem in its fundamental group may be found.  Indeed, the manifold 
may be first
decomposed into prime pieces (see for instance \cite{jactol:decomp}), whose fundamental
 groups form the free factors of
\pim, and each of these pieces either has automatic fundamental group or is
modeled on nilgeometry or solvgeometry (see \cite{ep:word}, Chapter 12 for a 
proof).  But any compact manifold modeled on nilgeometry or solvgeometry is 
finitely
covered by a circle bundle over a torus or a torus bundle over a circle,
respectively (see \cite{scott:geom}).  In particular, such a manifold is virtually
Haken, so Waldhausen's solution may be applied.  If the geometry of the
individual pieces is unknown at the outset, the search for an automatic
structure and the search for a finite cover which is Haken may be conducted in
parallel on each piece until one or the other is found.  A procedure to find the
automatic structure of a group is given in \cite{ep:word};  one to determine
whether a $3$--manifold is Haken is given in \cite{jactol:decomp}.

Although no known algorithm
 decides whether or not a $3$--manifold is geometric,  
some particular geometric manifolds can be recognized
algorithmically.  In \cite{thompson:s3} Thompson gives an algorithm to 
recognize the
$3$--sphere (see also \cite{rubinstein:s3}).
Similar methods provide recognition algorithms for the lens
spaces \cite{stocking:almnorm}.  Given the first homology of $M$
(which is easily computed) there is a finite list of lens spaces which $M$
could be homeomorphic to, so we can determine whether a manifold is a lens
space.  Rubinstein and Rannard have claimed to be able to recognize small
Seifert fibered spaces.  Note that the homeomorphism problem for Haken 
$3$-manifolds is decidable by work of Haken and Hemion
\cite{waldhausen:homeo},\cite{hemion:conjugacy}.  The
homeomorphism problem is also decidable for ``rigid weakly geometric''
$3$-manifolds by the work of Sela \cite{sela:isomorphism}.  An irreducible
$3$--manifold $M$ is \emph{rigid weakly geometric} if one of the following holds:
(1) $M$ is Haken,
(2) $M$ is geometric,
or (3) \pim is word hyperbolic and any irreducible $3$--manifold homotopy
equivalent to $M$ is in fact homeomorphic to $M$.

If a closed orientable $3$-manifold is known beforehand to satisfy the 
geometrization conjecture,
then the problem of determining its decomposition and the geometric structures
on the pieces is completely solved.  This is the case if, for instance, the
manifold is Haken (see \cite{thurston:bams}) or admits a symmetry with 
$\geq 1$-dimensional fixed point set (see \cite{chk:orbifold}).  For suppose
$M$ is such a manifold, which we may suppose is irreducible.  If $M$ is Haken,
there is an algorithm (algorithm 8.1 of \cite{jactol:decomp}) to find the
characteristic submanifold $\Sigma$.  
There are then three cases: (1) $\Sigma =
\emptyset$, (2) $\Sigma = M$ (in which case $M$ is a Seifert fibered space), 
and (3) $\Sigma$ is neither empty nor all of $M$.

In case (1) $M$ must be hyperbolic or a small Seifert fiber space.  
Sela gives 
a way to determine whether $M$ is a small
Seifert fiber space, and if so which one (details are in section 10 of
\cite{sela:isomorphism} and in Sela's Ph.D. thesis \cite{sela:thesis}).  
If $M$ is hyperbolic, the fundamental group is automatic, so a solution to the
word problem may be readily found.  The methods of the present paper may then
be used to find the hyperbolic structure on $M$.

In case (2), the algorithm from \cite{jactol:decomp} gives enough information to
determine the homeomorphism type of the Seifert fiber space, since it
constructs the fibering explicitly.  The geometric structure is then immediate.

In case (3), we first check whether $\partial \Sigma$ cuts $M$ into pieces
which are all homeomorphic to $T^2 \times I$ or an $I$-bundle over the Klein
bottle.  In this case $M$ has a solvgeometric structure which can be easily
determined.  Otherwise, $M\setminus \partial \Sigma$ consists of pieces which
are geometric or $I$--bundles which separate the geometric pieces from one
another.  Each geometric piece is either Seifert fibered or hyperbolic.
The Seifert fibered pieces
can be handled as in case (2), and the hyperbolic structures on the remaining
pieces can be computed by a variety of methods, for instance by a variation on
the methods of the present paper.  

Alternatively, results of and Petronio and Weeks
\cite{petronioweeks:flat} have shown that a manifold with torus boundary
components admits a finite volume hyperbolic structure if and only if Thurston's
hyperbolization equations for some ideal triangulation have a solution
corresponding to a combination of positively oriented and flat tetrahedra.
Any ideal triangulation is
accessible from any other by a sequence of finitely many moves, so this 
gives an algorithm
to find the hyperbolic structure, if one is known to exist.

We summarize the preceding discussion as:
\begin{theorem}
If a closed orientable 
$3$--manifold is known to satisfy the geometrization conjecture,
there are algorithms to find its canonical decomposition, and the geometric
structures on the pieces are algorithmically constructible.
\end{theorem}

\subsection{Outline}
We describe an algorithm to decide if a closed triangulated 
orientable $3$--man\-if\-old $M$ is
hyperbolic, assuming that the word problem in \pim\ can be solved.  
The key to the effectiveness of the algorithm is the following theorem of Gabai,
Meyerhoff and
N~Thurston.
\begin{theorem}\label{th:homhom}
{\rm\cite{ga:hom}}\qua Let $f\co M\to N$ be a homotopy equivalence between closed
$3$--manifolds,
where $M$ is irreducible and $N$ is hyperbolic.  Then $f$ is homotopic to a
homeomorphism.
\end{theorem}

This theorem implies that in order to determine that a closed $3$-manifold
$M$ 
is hyperbolic, we need only check (1) $M$ is irreducible, and (2) $\pim
\cong G$, where $G$ is a cocompact Kleinian group which acts freely on
hyperbolic space.  
For in this case $M$ and $\mathbb{H}^3/G$ are
Eilenberg--MacLane spaces with isomorphic fundamental groups, and therefore
homotopy equivalent.

\begin{lemma}\label{lemma:irreducible}
The irreducibility of $M$ is decidable.
\end{lemma}
\begin{proof}
This follows directly from classical work of Haken
and the recognizability of the $3$--sphere, as follows.  
The first step is to find a complete system of
$2$-spheres in $M$.  A \emph{complete system} is one which is guaranteed to
contain enough spheres to completely decompose $M$ into prime factors; 
the system
may contain redundant spheres.  The earliest algorithms to do this are due to
Haken;  more recent improvements can be found, for instance, in
\cite{jactol:decomp}.
For each sphere $S$ in the list, cut
$M$ along $S$.  If the resulting manifold is connected, then $S$ is a
non-separating sphere, and so $M$ is reducible.  If $M\setminus S$ is 
disconnected, let $M_1$ and $M_2$ be the pieces.  Each $M_i$ has a single
sphere boundary component; let $\hat{M}_i$ be the manifold
obtained by capping off this boundary component with a $3$-ball.  Now
Rubinstein and Thompson's $3$--sphere recognition algorithm may be applied to 
the $\hat{M}_i$, and
$S$ is a reducing sphere if and only if neither $\hat{M}_1$ or $\hat{M}_2$ is $S^3$.
\end{proof}

We can therefore assume for the remainder of the paper that $M$ is irreducible.
The question remains whether \pim\ has a discrete faithful
representation into $PSL_2\C$, with image which acts freely on hyperbolic space.  
In Section \ref{sec:alg} we give a procedure to find a finite list of candidate
representations, which  must contain a discrete faithful representation if
one exists.  We find this list by decomposing the representation variety into
irreducible components, and taking one representation from each component of the
appropriate dimension.  We also make some remarks on computations involving
algebraic numbers, which are used heavily in the sequel.

The next two sections are devoted to determining whether or not a given
candidate representation  gives rise to a hyperbolic structure on $M$.
The procedure given in Section \ref{sec:no} detects representations which
fail to be discrete and faithful with torsion--free image.  
The procedure given in
Section \ref{sec:yes} detects representations which are discrete and faithful
with torsion--free image, by
constructing a fundamental domain for the induced action on hyperbolic space.  
Finally in Section \ref{sec:structure} we prove the main theorems.  

\subsection{Acknowledgements}
Many thanks to Andrew Casson for coming up with the idea on which this paper is
based, and to my advisor Daryl Cooper for his enthusiasm and patience and for
suggesting that this be written up.  Thanks also to the NSF, which provided
partial support while this work was being done.

\section{\Qbar\ and some algebraic geometry in a nutshell}\label{sec:alg}

\subsection{Computations in \Qbar}
The algorithms given in this paper rely heavily on the ability to do
computations with algebraic numbers over the rationals.  Most of what we need to
use can be found in \cite{loos:ext}.
For the purposes of
computation, each algebraic number is to be thought of as a minimal polynomial
together with an isolating interval (or rectangle) with rational endpoints.  A
method for finding arbitrarily small isolating intervals for the solutions of a
polynomial over $\Q$ can be found for instance in 
Section 6.2 of \cite{bw:gro}.  The method
given there is for isolating real roots, but (as pointed out in an exercise) 
the method extends to complex roots as well.  
Since minimal polynomials are unique
up to multiplication by an element of $\Q$, this representation of algebraic
numbers gives us an effective test for equality.  

The reader who does the exercise from \cite{bw:gro} will
note that extracting the real or imaginary part of an algebraic number is then
straightforward.

For real algebraic numbers, we will need to decide the truth of statements of
the form ``$x<y$''.  Using the SIMPLE algorithm of \cite{loos:ext}, 
we may find a simple
extension $\Q(\alpha)$ of the rationals which contains both $x$ and $y$.  The
algorithm in Section 1.3 of \cite{loos:ext} may then be used to determine the
sign of $y-x$.    

We will also need to be able to add, multiply and take inverses and square roots
of algebraic numbers.  Algorithms to perform these operations are given in
\cite{loos:ext}.  The general idea is illustrated by the case of finding 
$\alpha \beta$, where $\alpha$ and $\beta$ are given.  First a polynomial is
found which must have $\alpha \beta$ as a root.  If $A(x)$ is the minimal
polynomial for $\alpha$ and $B(x)$ is the minimal polynomial for $\beta$, then
Loos shows that $\alpha \beta$ is a root of 
$r(x) = \mathrm{res}(y^m A(x/y), B(y))$, where $\mathrm{res}$ denotes the
resultant.  The irreducible factors of $r(x)$ are computed, and then the roots
of $r(x)$ are separated from one another by isolating intervals or rectangles.
Interval arithmetic is used to try to pick out the appropriate root of $r(x)$.
If the isolating intervals for $\alpha$ and $\beta$ are too big, then this can
fail on the first try.  In that case, the isolating intervals for $\alpha$ and
$\beta$ are shrunk until their product picks out a unique root of $r(x)$.

\subsection{Finding the candidate representations}\label{subsec:reps}
We now establish the following:
\begin{lemma}\label{lemma:reps}
There is an algorithm which takes as input a triangulation of a closed
$3$--manifold $M$, and outputs a finite list of rigid representations of 
\pim\ into $SL_2\Qbar$, where $\Qbar \subset \C$ is the algebraic closure of $\Q$.
Furthermore, if $M$ is hyperbolic, then this list contains a discrete faithful
representation.
\end{lemma}
\begin{proof}
From the triangulation of $M$ we obtain some presentation of $\pim$: 
\begin{equation}\label{eqn:pim}
\pim = \langle g_1,\ldots,g_n | w_1,\ldots,w_m \rangle
\end{equation}
  We then consider the representation
variety of homomorphisms 
\[R = \{\rho \co \pim  \to SL_2\C\} \subseteq \C^{4 n}\]
as defined for
example, in \cite{cs:var}.  Since \pim\ is finitely presented, this variety is
constructible as the zero set of a finite set of polynomials
with integer coefficients.  Indeed, $SL_2\C$ can be thought of as the variety 
in $\C^4$ defined by \mbox{$\{z_1 z_4 - z_2 z_3 = 1\}$}.  
The representation variety
$R$ is the subvariety of 
\[{SL_2\C}^n = \{(z_{1,1},\ldots,z_{4,1},\ldots,z_{1,n},\ldots,z_{4,n})\ |\ 
z_{1,i} z_{4,i} - z_{2,i} z_{3,i} = 1\ \forall i\} \]
obtained by adjoining the $4 m$
polynomials coming from the matrix equations corresponding to the 
relations of \pim.

\newcommand{\Rkl}{\ensuremath{R_{\{k,l\}}}}
For each pair $\{g_k,g_l\}$ of generators of \pim, we construct the subvariety
$\Rkl\!\subset R$ consisting of those representations for which the equations
\begin{equation}\label{eq:constraint}
\rho (g_k) = \left[ \begin{array}{cc} \ast & 1 \\ 0 & \ast
						\end{array} \right]
\mathrm{,\ and\ } 
\rho (g_l) = \left[ \begin{array}{cc} \ast & 0 \\ \ast & \ast
						\end{array} \right].
\end{equation}
hold.  In the notation introduced above, (\ref{eq:constraint}) 
corresponds to the equations 
$z_{2,k}=1$, $z_{3,k}=0$, and $z_{2,l}=0$.    Recall that a subgroup of $SL_2\C$ is said to be nonelementary if no
element of $\H^3\cup\partial \H^3$ is fixed by then entire subgroup.
Notice that if $\rho$ is any
representation so that $\langle \rho(g_k), \rho(g_l) \rangle
\subset SL_2\C$ is nonelementary, then $\rho$ is conjugate to a representation 
in $\Rkl$.  

Also observe that given some
representation $\rho$ of \pim\ which does not send $g_l$ to the identity, 
there are at
most four representations conjugate to $\rho$ which satisfy equations
(\ref{eq:constraint}).  Indeed, if $\rho$ and $C^{-1}\rho C$ both satisfy 
(\ref{eq:constraint}) for some $C\in SL_2\C$, then 
$C\left[\begin{array}{c} 1 \\ 0 \end{array} \right]$
must be an eigenvector of $\rho(g_k)$.
Likewise, $C\left[\begin{array}{c} 0 \\ 1 \end{array}\right]$
must be an eigenvector of $\rho(g_l)$.  Up to scaling the columns of $C$, 
there are at most four ways to arrange this.  Given one of the four choices, 
the requirements that the
determinant of $C$ be one and that $C^{-1}\rho(g_k)C$ have a one in the upper
right corner determine $C$ up to multiplication by $\pm I$. 

For each pair $\{g_k,g_l\}$, we would like to add the isolated points of 
\Rkl to our list of candidate representations.
Note that 
\Rkl is given by a set of polynomials with integer coefficients.  Let $I$ be
the ideal in $\Q[z_{1,1},\ldots,z_{4,n}]$ generated by these polynomials.
There exist algorithms (see, for instance \cite{bw:gro}, Theorem 8.101) 
to decompose this ideal into prime ideals.  Proposition 9.29 of
\cite{bw:gro} gives an algorithm to determine the dimension of such an ideal.
Each isolated point in \Rkl is part of the variety determined by some
zero--dimensional ideal $J$.  Let \mbox{$\{\zeta_{1,1},\ldots \zeta_{4,n}\}$}
be some such solution.  
Each coordinate $\zeta_{i,j}$ is a root of the monic polynomial $p(z_{i,j})$
which generates the
elimination ideal $J\cap \Q[z_{i,j}]$.  An algorithm to find this
polynomial is given in Section 6.2 of \cite{bw:gro} (see also Lemma 6.50 of
\cite{bw:gro}).  Possibilities for $\zeta_{i,j}$ are determined by finding the
irreducible factors of $p$ and isolating intervals or rectangles for the various
roots.  There is now a finite list of possible values for each $\zeta_{i,j}$.
Not every combination of choices necessarily corresponds to a representation, so
we must check the original relations in \pim\ by matrix computations over
$\Qbar$, for each of the finite number of combinations of choices.  Those
$4n$--tuples of choices which satisfy the relations are the candidate
representations coming from \Rkl.  Since there are finitely many pairs
$\{g_k,g_l\}$, the process described ends up with a finite list of
representations.

Suppose that $M$ is hyperbolic.  Then there is a
discrete faithful representation $\rho_0$ of \pim\ into $PSL_2 \C$ 
which gives the
hyperbolic structure on $M$.  This representation lifts to $SL_2 \C$
as proved in \cite{cs:var}.   We abuse notation slightly by continuing to refer
to the lifted representation as $\rho_0$. 
Note that the image $\Gamma = \rho_0(\pim)\subset SL_2\C$ cannot be cyclic, 
and that every nontrivial element is loxodromic.
Therefore, there must be some pair of generators
$\{g_k,g_l\}$ so that $\rho_0(g_k)$ and $\rho_0(g_l)$ are loxodromics with
distinct axes.  As $\Gamma$ is discrete,
$\langle \rho_0(g_k), \rho_0(g_l) \rangle
\subset SL_2\C$ is nonelementary (see \cite{beardon:gdg}, Theorem 5.1.2).
Therefore $\rho_0$ is conjugate to
a representation $\rho_0'$ in \Rkl.  
It is a well--known consequence of Mostow Rigidity 
that $\rho_0$ is unique up to conjugacy and orientation, and that nearby
representations are necessarily conjugate, so $\rho_0'$ is an
isolated point in \Rkl.  The discrete faithful representation $\rho_0'$ must
therefore appear in our list of candidate representations.
\end{proof}
\begin{remark}
If $R_{\{k,l\}}$contains a higher dimensional curve of representations at least one of 
which does not send $g_l$ to $I$, then it follows from \cite{cs:var} that $M$
is Haken, and the question of hyperbolicity could be answered by checking
whether the characteristic submanifold is empty, as pointed out in the
introduction.
However, the methods of the current paper still apply.
\end{remark}

\section{A procedure for rejecting representations}\label{sec:no}

Suppose $M$ is an orientable closed $3$--manifold, and $\rho\co\to
SL_2\Qbar$ is a representation as might be produced by Lemma \ref{lemma:reps}.
Suppose for the moment that $\rho$ is discrete
and faithful and that its image $\Gamma$ acts freely on hyperbolic space.  We 
immediately know a number of things about $\rho$ and its image.  For instance:
\begin{enumerate}
\item $\rho$ is irreducible (in particular, no geodesic is preserved by
$\rho(\pim)$).
\item $\rho$ has trivial kernel.
\item $\Gamma$ contains no nontrivial elliptics or parabolics.  Equivalently,
no nontrivial element of $\Gamma$ has trace in the interval $[-2,2]$.
\item Given a point in hyperbolic space, the subgroup of $\Gamma$ generated by 
elements which move that point a very small distance is abelian (in fact 
cyclic).
\end{enumerate}
Properties (1) and (3) follow immediately from the fact that $\H^3/\Gamma$ 
is a closed
$3$--manifold, and that no nontrivial element of $\Gamma$ has fixed
points.  If property (1) were to fail, then $\Gamma$ would be
elementary, and thus could not be cocompact.  If property (3) were to fail,
$\Gamma$ would either not act freely (if there were elliptics) or not be
cocompact (if there were parabolics).  For more details see \cite{beardon:gdg},
especially Chapters 4 and 5.
Property (2) is obvious.  Property (4) follows from and is made precise by 
the Margulis Lemma:
\begin{theorem}[Margulis Lemma]\label{th:marg}
There exists a small positive constant $\mu$ (eg, $\mu=\MARGEST$ works, see
\cite{gromov:mnpc}, page 107) such that the subgroup
$\Gamma_\mu(M,x)\subset \pi_1(M,x)$ generated by loops 
based at $x\in M$ of length at most $\mu$ is abelian, 
for any hyperbolic three--manifold
$M$, and any $x\in M$.  In particular, if $M$ is closed, then $\Gamma_\mu(M,x)
\cong \Z$ or is trivial.  
\end{theorem}
Furthermore, at least one of these four properties must fail if $\rho$ 
is not both discrete and faithful:
\begin{lemma}\label{lemma:indiscrete}
If $G < SL_2\C$ is indiscrete, consists only of loxodromics, 
and does not preserve any
geodesic, then given $x\in \H^3$  and $\epsilon > 0$ there are noncommuting 
matrices which move $x$ a distance less than $\epsilon$.  
\end{lemma}
\begin{proof}

For $A = \left[ \begin{array}{cc} a & b \\ c &d \end{array}\right]$ a
two--by--two complex matrix,
define:
\[\|A\|=   
(|a|^2+|b|^2+|c|^2+|d|^2)^{1/2}\]
The norm  $\|\cdot\|$ gives rise to a metric
on the space of two--by--two complex matrices defined by 
$\delta(A,B)=\|A-B\|$.  
The restriction of this metric to $SL_2\C$ gives the same topology as the
ordinary one (see \cite{beardon:gdg} for details).
Therefore the hypothesis that $G$ is indiscrete is equivalent to the 
existence of nontrivial elements $A$ of $G$ so that $\|A-I\|$ is 
arbitrarily small.
It is shown in \cite{beardon:gdg} that if $A\in SL_2\C$ acts on the upper half
space model of hyperbolic
space in the usual way, and $j$ is the point
$(0,0,1)$, then $\|A\|^2 = 2 \cosh d(j,Aj)$, where $d$ is the hyperbolic
metric.  A simple calculation shows that if there are nontrivial matrices 
$A\in G$ with $\|A-I\|$ arbitrarily small, then there must also be matrices
with $\|A\|^2-2$ (equivalently $\cosh^{-1}(\frac{\|A\|^2}{2})$)
arbitrarily small.
Since the only elements of $SL_2\C$ which fix $j$ (and thus have $\|A\|^2-2=0$) 
are elliptics whose axes pass through $j$, this means there must be nontrivial
matrices in $G$ with $\cosh^{-1}(\frac{\|A\|^2}{2})$ arbitrarily small 
but positive.  We claim that we can find such matrices which do not
commute.  

Suppose we cannot.  Then all matrices in $G$ which move $j$ a sufficiently
small distance commute, and thus share a common axis, which we'll call
$\gamma$.  Let $G_{\gamma}$ be the subgroup of $G$ which sends $\gamma$ to itself.

Rotation angle and translation length are continuous functions on $SL_2\C$, so
we can find elements of $G_{\gamma}$ which both translate $\gamma$ an
arbitrarily small distance along itself and rotate $\H^3$ about $\gamma$ by an
arbitrarily small angle.

  By hypothesis, $G\setminus G_\gamma$ is
nonempty. Let $C\in G\setminus G_\gamma$.  The subgroup 
$ C (G_\gamma) C^{-1}$ of $G$ preserves the geodesic $C(\gamma) \neq \gamma$
and so its elements do not commute with those of $G_\gamma$.  Since rotation
angle and translation length are invariants of conjugacy class, there are
elements of $ C (G_\gamma) C^{-1}$ with arbitrarily small such.  There are
therefore elements of $ C (G_\gamma) C^{-1}$ which move $j$ as small a distance
as we like.  This contradicts our assumption that all matrices in $G$ which
move $j$ a sufficiently small distance commute.  

Finally let $x$ be an arbitrary point in $H^3$.  If $g$ is an element of
$SL_2\C$ which takes $x$ to $j$ we may apply the above argument to $g^{-1}Gg$.
\end{proof}

\begin{figure}[ht!]
\begin{center}
\begin{picture}(0,0)%
\includegraphics{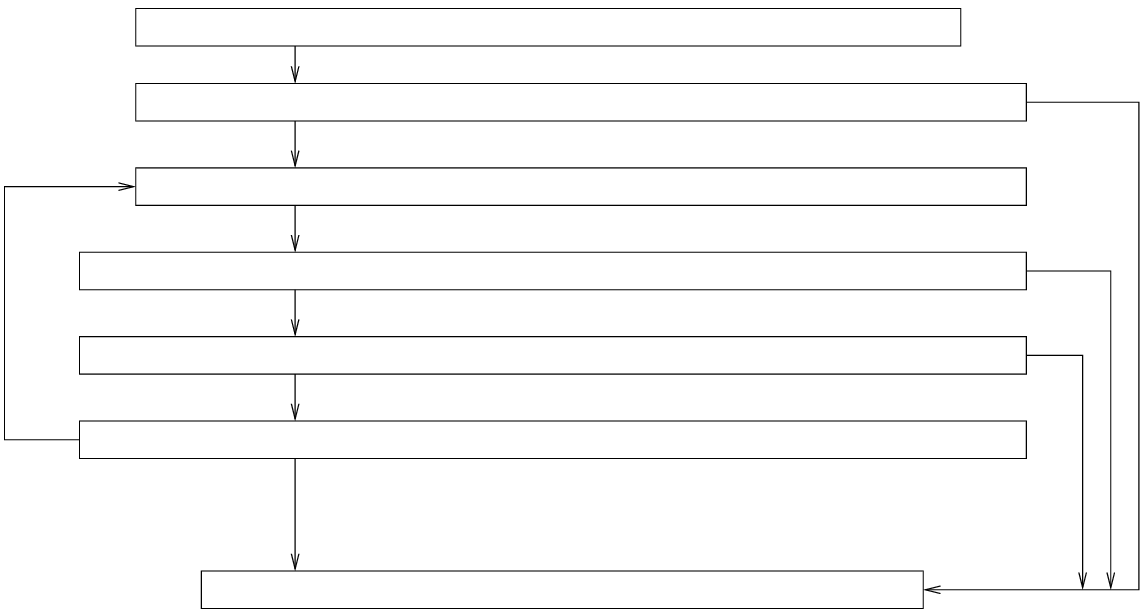}%
\end{picture}%
\setlength{\unitlength}{2368sp}%
\begingroup\makeatletter\ifx\SetFigFont\undefined%
\gdef\SetFigFont#1#2#3#4#5{%
  \reset@font\fontsize{#1}{#2pt}%
  \fontfamily{#3}\fontseries{#4}\fontshape{#5}%
  \selectfont}%
\fi\endgroup%
\begin{picture}(9099,4824)(1414,-5923)
\put(2551,-1936){\makebox(0,0)[lb]{\smash{\SetFigFont{7}{8.4}{\familydefault}{\mddefault}{\updefault}
Is $\rho$ reducible?(Lemma \ref{lemma:reducible})
}}}
\put(2551,-1336){\makebox(0,0)[lb]{\smash{\SetFigFont{7}{8.4}{\familydefault}{\mddefault}{\updefault}
Given some representation $\rho:\pim\longrightarrow SL_2\Qbar$
}}}
\put(2551,-2611){\makebox(0,0)[lb]{\smash{\SetFigFont{7}{8.4}{\familydefault}{\mddefault}{\updefault}
List some (more) elements of the free group on $G$.  Call the list $W$.
}}}
\put(3826,-2236){\makebox(0,0)[lb]{\smash{\SetFigFont{7}{8.4}{\familydefault}{\mddefault}{\updefault}
NO
}}}
\put(3826,-3586){\makebox(0,0)[lb]{\smash{\SetFigFont{7}{8.4}{\familydefault}{\mddefault}{\updefault}
NO
}}}
\put(3826,-4336){\makebox(0,0)[lb]{\smash{\SetFigFont{7}{8.4}{\familydefault}{\mddefault}{\updefault}
NO
}}}
\put(3826,-5011){\makebox(0,0)[lb]{\smash{\SetFigFont{7}{8.4}{\familydefault}{\mddefault}{\updefault}
NO
}}}
\put(9676,-3811){\makebox(0,0)[lb]{\smash{\SetFigFont{7}{8.4}{\familydefault}{\mddefault}{\updefault}
YES
}}}
\put(9676,-3136){\makebox(0,0)[lb]{\smash{\SetFigFont{7}{8.4}{\familydefault}{\mddefault}{\updefault}
YES
}}}
\put(9676,-1711){\makebox(0,0)[lb]{\smash{\SetFigFont{7}{8.4}{\familydefault}{\mddefault}{\updefault}
YES
}}}
\put(2101,-4636){\makebox(0,0)[lb]{\smash{\SetFigFont{7}{8.4}{\familydefault}{\mddefault}{\updefault}
Is every $w \in W$ so that $\rho(w)=I$ equivalent to $1$ in \pim? (word problem)
}}}
\put(2101,-3961){\makebox(0,0)[lb]{\smash{\SetFigFont{7}{8.4}{\familydefault}{\mddefault}{\updefault}
Are there small noncommuting matrices in $\rho(W)$? (Lemma \ref{lemma:noncommuting})
}}}
\put(2101,-3286){\makebox(0,0)[lb]{\smash{\SetFigFont{7}{8.4}{\familydefault}{\mddefault}{\updefault}
Does $\rho(W)$ contain nontrivial elliptics or parabolics? (Lemma \ref{lemma:trace})
}}}
\put(1501,-4486){\makebox(0,0)[lb]{\smash{\SetFigFont{7}{8.4}{\familydefault}{\mddefault}{\updefault}
YES
}}}
\put(3081,-5836){\makebox(0,0)[lb]{\smash{\SetFigFont{8}{9.6}{\familydefault}{\mddefault}{\updefault}
$\rho$ is not discrete and faithful with torsion--free image
}}}
\end{picture}
 
\caption{A procedure to reject representations.  $G$ is a finite generating set
for \pim.}
\label{figure:noflow}
\end{center}
\end{figure}

The procedure for rejecting representations, outlined in Figure
\ref{figure:noflow}, is based on the conditions (1)--(4), so it 
requires that we
be able to recognize when these conditions have failed.

\begin{lemma}\label{lemma:reducible}
There is an algorithm to decide whether or not a finite set of matrices in 
$SL_2 \Qbar$ share a common eigenvector.
\end{lemma}
\begin{proof}
Suppose 
$A=\left[ \begin{array}{cc} a & b \\ c &d \end{array}\right] \in
SL_2\Qbar\setminus \{\pm I\}$. 

If $c\neq 0$ then $A$ has
$\left[ \begin{array}{c} (a\pm \sqrt{(a+d)^2-4})/(2 c) \\ 
1 \end{array}\right]$ as eigenvectors.
If $c = 0$ then $A$ has eigenvectors 
$\left[ \begin{array}{c} 1 \\ 0 \end{array}\right]$ and
$\left[ \begin{array}{c} b/(a-d) \\ 1 \end{array} \right]$, unless $a-d=0$
in which case $A$ has 
$\left[ \begin{array}{c} 1 \\ 0 \end{array}\right]$ as its sole eigenvector.
Of course any multiple of an eigenvector is an eigenvector, but we have chosen
specific ones precisely to remove that ambiguity.  Since we can do arithmetic
and square roots in $\Qbar$ and compare elements of $\Qbar$, it is possible to 
decide whether or not two vectors of the given form are the same.
\end{proof}

\begin{lemma}\label{lemma:trace}
There is an algorithm to decide whether or not a finite set of matrices $S \in
SL_2\Qbar$ contains any nontrivial elliptics or parabolics.
\end{lemma}
\begin{proof}
 A matrix in $SL_2\C$ acts elliptically or parabolically on hyperbolic space if
and only if its trace is in the interval $[-2,2]$. 
For each $A\in S$ we can find the algebraic number $\mathrm{trace}(A)$.  As
discussed in Section \ref{sec:alg}, there are algorithms to find the imaginary
part of $\mathrm{trace}(A)$, and to test equalities and inequalities, so it can
be determined whether or not 
$\mathrm{trace}(A) \in [-2,2]$.  If $\mathrm{trace}(A) =  2$, 
we may then check whether or not $A = I$.  
\end{proof}

To apply the discreteness criterion effectively, we need a checkable condition
on the elements of $\Gamma$.  As noted in the proof of Lemma
\ref{lemma:indiscrete}, a matrix $A$ acting on hyperbolic
space moves $j$ a hyperbolic distance of $\cosh^{-1}(\frac{(\|A\|)^2}{2})$.
If $\Gamma$ is discrete, then the Margulis Lemma implies that all matrices $A$
with $\cosh^{-1}(\frac{(\|A\|)^2}{2})<\MARGEST$ must commute with each other.
To avoid having to calculate with transcendental functions and numbers, we note
that, for instance $\|A\|^2<2+2^{-58}$ implies 
$\cosh^{-1}(\frac{(\|A\|)^2}{2})<\MARGEST$.

\begin{lemma}\label{lemma:noncommuting}
There is an algorithm to decide whether or not a finite set of matrices 
$S\in SL_2\Qbar$
contains noncommuting elements $A$ and $B$ so that  $\|A\|-2 <
2^{-58}$ and $\|B\|-2 < 2^{-58}$.
\end{lemma}
\begin{proof}
For each $A\in S$, we compute $\|A\|$ (which involves just multiplication,
addition, and a square root).  For each pair of matrices $A$, $B$ of $S$ so 
that $\|A\|-2 <2^{-58}$ and $\|B\|-2 < 2^{-58}$ we check 
whether or not $AB-BA$ is the zero matrix.
\end{proof}

Now we can prove the main result of this section.
\begin{theorem}\label{th:no}
There exists an algorithm which takes as input a 
triangulated orientable closed 
three-manifold $M$, a
solution to the word problem in $\pim$ and a representation $\rho\co\pim\to
SL_2\Qbar$, and decides if the representation is \emph{not} discrete and
faithful with torsion--free image.
\end{theorem}
\begin{proof}
The idea is to detect the failure of conditions (1)--(4) listed at the beginning
of this section.  By Lemma \ref{lemma:indiscrete}, one of these conditions must
fail if $\rho$ is not discrete and faithful with torsion--free image.

We start with condition (1).   Let $G=\{g_1,\ldots,g_n\}$ from the
presentation (\ref{eqn:pim}) in Section \ref{subsec:reps}.
The representation $\rho$ is reducible if and only if all the elements of
$\rho(G)$
share a common eigenvector.  This condition is checkable by Lemma
\ref{lemma:reducible}.  

If the representation is irreducible, we begin systematically listing elements
of the free group $\mathcal{F}(G)$ in such a way that any element is eventually
listed.  For each finite list produced this way, we look for evidence that
conditions (2)--(4) have failed for $\rho$.  Let $W$ be the finite list at some
stage in the listing process.  By performing matrix multiplications over 
$\Qbar$, we can determine $\rho(w)$ for any $w\in W$.  Let $P = \rho(W)$.

If condition (2) fails, then there will eventually be elements of $W$ which are
not equivalent to the identity in \pim, but which are in the kernel.
If $\rho(w)=I$ for some $w \in W$, then we apply the solution of the word 
problem in \pim\ to
determine whether or not $[w]$ is the identity element in \pim.  
Thus if condition (2) fails, we will eventually discover it.
Note that this is the
first use of the assumption that \pim\ has solvable word problem.

If condition (3) fails, then eventually $P$ will contain nontrivial elliptics or
parabolics.  By Lemma \ref{lemma:trace}, we can detect this.

If condition (4) fails, then there are
arbitrarily small noncommuting matrices.  In particular, $P$ will eventually
contain noncommuting matrices $A$ and $B$ with $\|A\|-2 <2^{-58}$ and 
$\|B\|-2 < 2^{-58}$.  By Lemma \ref{lemma:noncommuting}, we can detect such
matrices.

If $\rho$ is discrete and faithful with torsion--free image, then none of the
conditions (1)--(4) can fail, and so this procedure will continue endlessly.
Otherwise, one of the four conditions will eventually be discovered to fail and
the procedure will stop.
\end{proof}

\section{Recognizing a discrete faithful representation}\label{sec:yes}

In this section we give a procedure which takes as input a representation
\mbox{$\rho\co\pim\to SL_2\Qbar$} and decides if
$\rho$ is discrete and faithful with torsion free image.  The plan of attack is to first try to build a 
fundamental region for the action,
\`{a} la Riley in \cite{riley:p}, and then to use the solution to the word
problem to verify that $\rho$ is injective.  

It is convenient to consider both the ball and upper half--space models of
hyperbolic space as subsets of the quaternions (this is a point of view taken in
\cite{beardon:gdg} and in \cite{riley:p}).  Specifically we can declare
\mbox{$\B^3=\{a+bi+cj\,|\, a^2+b^2+c^2<1\}$} and 
\mbox{$\H^3=\{a+bi+cj\,|\,c>0\}$} and
agree to always pass between them by way of the isometry 
$f\co\H^3\to\B^3$ given by $f(w)=(w-j)(w+j)^{-1}j$.
A matrix 
$A=\left[ \begin{array}{cc} a & b \\ c &d \end{array}\right] \in SL_2\C$ 
acts on $\H^3$ as $q \mapsto (aq + b)(cq +d)^{-1}$, and acts on $\B^3$ by 
$f \circ A \circ f^{-1}$.  The action on $\B^3$ extends to a M\"obius
transformation of $\R^3 \cup \{\infty\}$, where $\R^3 = \{a+bi+cj+dk|d=0\}$, and 
if $A\in SL_2 \Qbar$, then this action preserves the set $\Qbar^3\cup
\{\infty\}$.  

So long as $A$ is not elliptic or $\pm I$ there is a unique Euclidean sphere
in $\R^3$  on
which the M\"obius transformation acts isometrically (ie, the determinant of
the Jacobian matrix is 1).  We call this sphere the
\emph{isometric sphere of $A$}.  
\begin{definition}
For a point $v\in\R^3$ and a sphere $S$ with center $c$ and radius $r$ we say
$v$ is \emph{inside} $S$ if $\|v-c\|<r$, $v$ is \emph{outside}
 $S$ if $\|v-c\|>r$, and $v$ is \emph{on}
$S$ if $\|v-c\|=r$.  Note that if $r$ and the coordinates of the points $c$ and
$v$ are all algebraic, then each of the three conditions  is   
 decidable.  
\end{definition}
Note that since the intersection of an isometric sphere with $\B^3$ 
is a hyperbolic plane in
$\B^3$, the intersection of the interior or exterior of an isometric sphere with
$\B^3$ is convex with respect to the hyperbolic metric.

The isometric spheres are not intrinsic to the hyperbolic metric in the sense 
that, say, the axis of a
loxodromic is intrinsic.  However, they give 
us a convenient way to construct a particular fundamental domain for a discrete
action.  In the case that no group element acts as a nontrivial elliptic, the
closure of the set of points in $\B^3$ which are outside every isometric sphere
is a fundamental domain, sometimes called the \emph{Ford domain} (though this
term is sometimes used differently).  
In fact, this set is precisely the Dirichlet domain centered at 0 (see
\cite{mt:hmkg}, page 45).  
This is the fundamental domain which we will
construct for a discrete action.

According to Proposition 1.3 of \cite{mt:hmkg}, the isometric
sphere of an isometry $\phi\co\B^3\!\to\B^3$ has Euclidean center 
$\frac{\phi^{-1}(0)}{|\phi^{-1}(0)|^2}$ and Euclidean radius
$(\frac{1}{\phi^{-1}(0)^2}-1)^{1/2}$.  In particular, we have:
\begin{lemma}\label{lemma:spherecoords}
If $S_A$ is the isometric sphere of $A\in SL_2\Qbar$ acting on $\B^3$, and $S_A$
has radius $r_A$ and center $c_A = (x_A,y_A,z_A)$, then $r_A$,$x_A$,$y_A$ and
$z_A$ are algebraic numbers which can be computed from the entries of the matrix
$A$.
\end{lemma}

Let $W$ be some finite list of words in $\mathcal{F}(G)$, closed under inverses,
and let $P = \rho(W)\setminus\{\pm I\}$.  
We say that $P$ \emph{satisfies condition $NE$} if no
nontrivial element of $P$ is an elliptic, and no two nontrivial elements of $P$
have the same isometric sphere.
\begin{lemma}\label{lemma:NE}
There is an algorithm to decide whether or not $P$ satisfies condition $NE$.
\end{lemma}
\begin{proof}
This follows from our ability to compute traces and from Lemma
\ref{lemma:spherecoords}. 
\end{proof}

\begin{figure}[ht!]
\begin{center}
\begin{picture}(0,0)%
\includegraphics{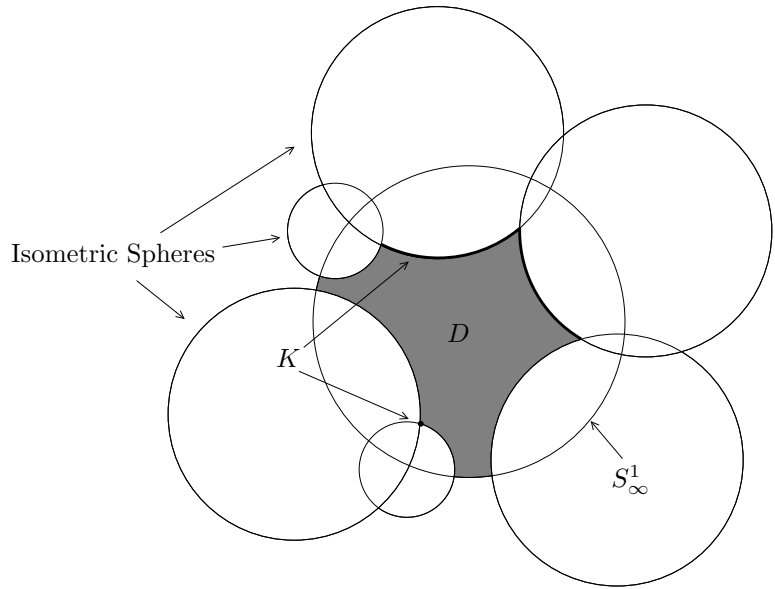}%
\end{picture}%
\setlength{\unitlength}{1579sp}%
\begingroup\makeatletter\ifx\SetFigFont\undefined%
\gdef\SetFigFont#1#2#3#4#5{%
  \reset@font\fontsize{#1}{#2pt}%
  \fontfamily{#3}\fontseries{#4}\fontshape{#5}%
  \selectfont}%
\fi\endgroup%
\begin{picture}(11933,9119)(226,-9068)
\put(9611,-7494){\makebox(0,0)[lb]{\smash{\SetFigFont{10}{12.0}{\familydefault}{\mddefault}{\updefault}
\small $S^{1}_\infty$
}}}
\put(7031,-5214){\makebox(0,0)[lb]{\smash{\SetFigFont{10}{12}{\familydefault}{\mddefault}{\updefault}
\small $D$
}}}
\put(4351,-5611){\makebox(0,0)[lb]{\smash{\SetFigFont{10}{12}{\familydefault}{\mddefault}{\updefault}
\small $K$
}}}
\put(226,-3961){\makebox(0,0)[lb]{\smash{\SetFigFont{10}{12.0}{\familydefault}{\mddefault}{\updefault}
\small Isometric Spheres
}}}
\end{picture}
\caption{This is what $K$ and $D$ might look like in dimension 2.  In this
example $K$ has four vertices and two edges.}
\label{figure:kpic}
\end{center}
\end{figure}

If $P$ satisfies condition $NE$, let $D$ be the closure in $\B^3$ of the 
intersection of the exteriors of the
isometric spheres $S_A$ for $A\in P$.  Then $D$ is convex with respect to the
hyperbolic metric, and 
bounded by finitely many (possibly noncompact or infinite area) hyperbolic
polygons, each of which is a subset of some isometric sphere.  Define a
\emph{vertex} of $D$ to be a point in $D$ which is contained in at least three
distinct isometric spheres.
Let $K$ be the
intersection of the convex hull of the  vertices of $D$ with
$\partial D$. Figure \ref{figure:kpic} gives an
idea what $K$ and $D$ might look like one dimension down (in dimension 2
vertices need only be contained in two distinct isometric spheres).  
$K$ can be thought of
as the ``compact part'' of the boundary of $D$. In particular, if $D$ is
compact, then $K = \partial D$. 
$K$ has an obvious cell decomposition.  In particular, the vertices of $K$ are
precisely the vertices of $D$.  If a pair of vertices are contained in two or
more distinct isometric spheres, then their convex hull forms an edge of $K$.
The $2$--cells of $K$ are the convex hulls of maximal collections of at least
three vertices, subject to the condition that the convex hull is contained in
$K$.  
Every cell in $K$ is the convex hull of its vertices, and is thus completely
determined by them.
\begin{lemma}\label{lemma:K}
If $P$ satisfies condition $NE$, there is an algorithm to construct $K$.  By
constructing $K$ we mean computing the locations of the vertices of $K$, and
specifying the edges and faces of $K$ as subsets of the vertex set of $K$.
\end{lemma}
\begin{proof}
Given three spheres in Euclidean space, it is not hard to write down a set
of inequalities in terms of their centers and radii which are satisfied if and
only if the spheres intersect in a pair of distinct points.  Moreover, these
inequalities involve only arithmetic and square roots, and so may be checked
algorithmically if the centers and radii are given by algebraic numbers.  If the
spheres intersect in a pair of points, then the coordinates of these points can
be explicitly written down (again using only arithmetic and square roots) in
terms of the centers and radii of the spheres.  

Therefore, for each triple of
isometric spheres of elements of $P$
it is possible to decide if they intersect in $\B^3\subset \R^3$, and if so to
write down a triple of algebraic numbers which are the coordinates of their
intersection.  These triples are potential vertices of $K$.  

Let $\mathcal{V}$ be
the set of potential vertices of $K$.   Since every vertex of $K$ is contained
in at least three distinct isometric spheres, every vertex of $K$ is in
$\mathcal{V}$.
However $\mathcal{V}$ may contain vertices not in $D$.  For each $A\in P$,
$v\in\mathcal{V}$, we check that $v$ is outside or on the sphere $S_A$.  
If $v$ is
inside $S_A$, $v$ is outside $D$, and so it is
discarded.  After discarding these vertices, we continue to refer to the 
smaller  set of vertices as $\mathcal{V}$.  After all sphere--vertex pairs have
been checked, this set is equal to the 
zero--skeleton of $K$. 

Now we find the edges of $K$ by looking for pairs of spheres whose intersection
contains two vertices. 
For each pair of vertices $\{v,w\}$ in $\mathcal{V}$, we can construct the list
of matrices $A$ in $P$ so that $v$ and $w$ are both on $S_A$.  There is
an edge connecting $v$ to $w$ in $K$ if
and only if this list contains at least two elements.  
Let $\mathcal{E}$ be the set of unordered pairs of vertices contained in at 
least two
isometric spheres.  Clearly every edge of $K$ is
represented by some element of $\mathcal{E}$.  Conversely, suppose
$\{v,w\}\subset S_A\cap S_B$ for some $A$,$B\in P$.  
Since $D$ is convex, the geodesic segment
between $v$ and $w$ is  contained in $D$.  Since $S_A\cap S_B$ is convex,
this segment is in fact contained in the boundary of $D$, and is therefore part
of $K$.  Thus we can construct the $1$--skeleton of $K$.

It remains only to find the $2$--cells of $K$.
For each $A\in P$, the set $\mathcal{V}\cap S_A$ can be constructed.  
There is a
$2$--cell  of $K$ formed from a part of $S_A$ if and only if the set 
$\mathcal{V}\cap S_A$ contains
at least three vertices and these vertices are linked up by edges in
$\mathcal{E}$ to form a circuit.  These conditions can be checked on $S_A$ for
each $A\in P$, so we can produce a list of faces
$\mathcal{F}$.  This completes the combinatorial data needed to construct $K$. 
\end{proof}
\begin{remark}(Notation)\qua
It is no extra work in the above to keep track of which matrix corresponds to
which face.  Since each isometric sphere contributes at most one face to
$\mathcal{F}$ and since $P$ satisfies condition $NE$, we may unambiguously refer
to a face as $F_A$ if it is contained in the isometric sphere of $A$.
\end{remark}
\begin{lemma}\label{lemma:twosphere}
There is an algorithm to decide whether or not $K$ is a $2$--sphere.
\end{lemma}
\begin{proof}
To check that $K$ is connected, it suffices to check that any pair of vertices
in $\mathcal{V}$ can be connected by a sequence of edges.  
The Euler characteristic of $K$ is equal to 
$|\mathcal{F}|-|\mathcal{E}|+|\mathcal{V}|$. 
\end{proof}

Note that if $K$ is a $2$--sphere, then $D$ is a compact finite sided 
hyperbolic polyhedron with $K=\partial D$.
If $D$ is a fundamental domain for a discrete action of some group on 
hyperbolic space, then we should have $A(F_A)=F_{A^{-1}}$ for every 
$F_A \in \mathcal{F}$.  If this is the case, then $K$ can be given the
additional structure of a polyhedron with face identifications.
\begin{lemma}\label{lemma:identifiable}
If the complex $K$ produced by Lemma \ref{lemma:K} is a $2$--sphere, then there
is an algorithm to
check whether or not $A(F_A)=F_{A^{-1}}$ for every  $F_A \in \mathcal{F}$.
\end{lemma}
\begin{proof}
We first check, for each $F_A \in \mathcal{F}$, that we also have
$F_{A^{-1}}\in\mathcal{F}$.
Since faces are uniquely determined by their vertices, it then suffices to
check that $A$ sends the vertices of $F_A$ to those of $F_{A^{-1}}$.  Since we
have kept track of the coordinates of the vertices, and can explicitly 
compute the action of $A$ on points in $\B^3$ with algebraic coordinates, 
this is straightforward.
\end{proof}

The question of when a polyhedron with face identifications is the fundamental
domain of a discrete action consistent with those face identifications is
answered by the Poincar\'e Polyhedron Theorem. 
We use the theorem in the following specialized form (more general versions
can be found in the standard references \cite{maskit:p} and \cite{seifert:p}):
\begin{theorem}[Poincar\'e Polyhedron Theorem]\label{th:ppt}
Let $D$ be a compact polyhedron in hyperbolic space,
together with face pairings which are orientation-preserving isometries
generating the group $\Gamma$.  Suppose the angle sum about every image of an 
edge
in the identified polyhedron is $2\pi$. Then:
\begin{enumerate}
\item $D$ is a fundamental polyhedron for the (fixed--point free)
action of $\Gamma$, and
\item $\Gamma$ has a presentation as an abstract group with the generators
equal to the face pairings and the relations precisely the edge cycles.
\end{enumerate}
\end{theorem}
\begin{remark}
Theorem \ref{th:ppt} gives a sufficient condition for a compact convex
polyhedron with face pairings to
be the 
fundamental domain for a discrete cocompact action.  If we also require that the
action be fixed--point free (so that the quotient is a hyperbolic
$3$--manifold), then the sufficient condition is also clearly necessary.
\end{remark}
\begin{lemma}\label{lemma:fundamental}
If the complex $K$ produced by Lemma \ref{lemma:K} is a $2$--sphere which
satisfies $A(F_A)=F_{A^{-1}}$ for every  $F_A \in \mathcal{F}$, then there is an
algorithm to decide whether or not $K$ is the boundary of a fundamental domain 
for a discrete fixed--point free action on hyperbolic space, generated by
$\{A\in P\ |\ F_A \in \mathcal{F}\}$.
\end{lemma}
\begin{proof}
We use a combination of exact and numerical computations to check whether the
conditions of Theorem \ref{th:ppt} are satisfied at each edge.
Specifically, let $e_1$ be some edge of $K$.  Adjacent to $e_1$ are 
precisely two faces, $F_A$ and $F_B$.
Let $F_{A_1}=F_A$ be one of the two faces.
Let $e_2$ = $A_1(e_1)$.  By hypothesis, $e_2$ is an edge of $K$.  
The edge $e_2$
is adjacent to two faces, one of which is $A_1(F_{A_1})$.    
The other face is $F_{A_2}$ for some $A_2\in P$.
We inductively define $e_{i+1}$ to be $A_{i}(F_{A_i})$ and $F_{A_{i+1}}$ 
to be the face
adjacent to $A_i(F_{A_i})$ across $e_{i+1}$.  
Let $k$ be the first integer for which $e_k=e_1$, the edge we 
started with, in which case $A_{k-1}(F_{A_{k-1}})=F_B$ and $F_{A_k}=F_B$.

If $\{e_1,\ldots,e_{k-1}\}$ satisfies the edge
condition, then, in particular, $A_{k-1} \circ \cdots \circ A_1 = I$.  If the
exact calculation shows that this is the case, then the angle sum is  
some integer multiple of $2\pi$.  

Suppose $S_1$ and $S_2$ are two spheres whose intersection contains an edge.
Then
their dihedral angle is given by 
$\cos^{-1}(\frac{\|c_1-c_2\|^2-{r_1}^2-{r_2}^2}{2r_1 r_2})\in(0,\pi)$, if 
$S_i$ is the sphere of
radius $r_i$ centered at $c_i$.  In the case where $r_i$ and the coordinates of
$c_i$ are given by algebraic numbers (or any numbers which can be determined to
any desired accuracy) it is not hard to show
that the angle can then be determined numerically to any desired accuracy.
We can therefore do a
calculation of the angle sum around an edge with tolerance less than, 
say $\pi/2$.  
Together with the exact
calculation, this tells us whether the edge condition is satisfied.  
As there are a finite number of edges in the complex $K$, it will take a finite
amount of time to check them all.  
\end{proof}
\begin{remark}
The proof of the preceding lemma can be modified to deal with the case of an
arbitrary polyhedron with identifications, provided that all the identifying
isometries are given as matrices in $SL_2\Qbar$, and all vertices lie inside
$S_{\infty}^{2}$ and have algebraic coordinates.  The version given is sufficient
for our purposes, however.
\end{remark}
If $K$ as constructed in Lemma \ref{lemma:K}
is the boundary of a fundamental domain, we will refer to the group
generated by $\{A\in P\ |\ F_A \subset K\}$ as $\Gamma_K$. 
$\Gamma_K$ is clearly a subgroup of $\rho(\pim)$, and the quotient of
hyperbolic space by the action of $\Gamma_K$ is a closed hyperbolic
$3$--manifold.  
\begin{remark}
It is possible that $\Gamma_K$ is a proper subgroup of 
$\rho(\pim)$.  In the case that $\rho(\pim)$ is discrete, then the quotient of
hyperbolic space by $\Gamma_K$ is a finite cover of the quotient of hyperbolic
space by $\rho(\pim)$.
\end{remark}
\begin{lemma}\label{lemma:gammakisgamma}
If $K$, constructed as above from a representation $\rho\co\pim\to
SL_2\Qbar$, 
is the boundary of a fundamental domain $D$
of a discrete fixed--point free action
on hyperbolic space, then there is an
algorithm to decide whether or not
$\Gamma_K = \rho(\pim)$, where $\Gamma_K$ is the group
generated by the elements of $\rho(\pim)$ which pair the faces of $K$.
\end{lemma}
\newcommand{\zo}{\ensuremath{\mathbf{0}}}
\begin{proof}
For each generator $g$ of \pim, we wish to check whether $\rho(g)$ is in
$\Gamma_K$.  Let $\{T_1,\ldots,T_s\}$ be the generators of $\Gamma_K$, which are
also the face pairings of $D$.  Each face of $K$ is part of the
isometric sphere $S_T$ for some
$T\in\{T_1,\ldots,T_s\}\cup\{T_1^{-1},\ldots,T_s^{-1}\}$.

Note that $\zo=(0,0,0)$ is in the interior of $D$.
Since hyperbolic space is tiled by copies of $D$,
$\rho(g)(\zo)$ is in $\gamma(D)$ for some $\gamma \in \Gamma_K$.  
But the following are clearly equivalent:
\begin{enumerate}
\item $\rho(g)(\zo)\in \gamma(D)$
\item $\gamma^{-1}\rho(g)(\zo)\in D$ 
\item $\gamma^{-1}\rho(g)(\zo)$ is outside or on $S_{T}$ for 
all $T \in \{T_1,\ldots,T_s\}\cup\{T_1^{-1},\ldots,T_s^{-1}\}$.
\end{enumerate}
From our earlier remarks, this last statement is clearly checkable.
If $\rho(g) \in \Gamma_K$ then $\rho(g)=\gamma$ whenever the statement is
true. 

We can systematically list all possible words in the
free group $\mathcal{F}(T_1,\ldots,T_s)$ and for each word $w$ in the list 
check to see whether condition (3) above holds for
$w^{-1}\rho(g)(\zo)$.  
We must eventually find
such a $w$, since the copies of $D$ tile hyperbolic space.  When this $w$ is
found, test whether $w^{-1}\rho(g) = I$.  If so, then $\rho(g)\in \Gamma_K$, 
and we go on to the next generator.  
If every generator has image in $\Gamma_K$, then
$\Gamma_K = \rho(\pim)$.  Otherwise, $\Gamma_K$ is a proper subgroup of
$\rho(\pim)$.
\end{proof}
\begin{lemma}\label{lemma:faithful}
Suppose $\rho\co\pim\to SL_2\Qbar$, and $K$ is a complex constructed 
as above from a finite subset of $\rho(\pim)$, and that
$\Gamma_K= \pim$.  If there
is an algorithm to solve the word problem in \pim, then there is an algorithm to
decide whether or not $\rho$ is faithful.
\end{lemma}
\begin{proof}
We can now think of $\rho$ as an epimorphism of abstract groups 
\[\rho\co\langle g_1,\ldots,g_n|r_1,\ldots,r_m \rangle \longrightarrow
\langle T_1,\ldots,T_s|c_1,\ldots,c_r\rangle\] where the $c_i$ are words coming from the
edge cycles.  
Recall that we started with some finite set $W$ of words in the free group 
$\langle g_1, \ldots, \ldots, g_n\rangle$.  The $T_i$ correspond to the
matrices in $P=\rho(W)$ which pair the faces of the polyhedron $K$, and generate
$\Gamma_K$.
To check whether or not $\rho$ is faithful, we attempt to construct an inverse 
map.  For each $T_i$ pick some
$w_i$ in the finite set $\rho^{-1}(T_i)\cap W$.  
The $w_i$ thus chosen determine a homomorphism from the free group
\mbox{$\langle T_1,\ldots,T_s\rangle$}
 to \pim.  We can check whether this map factors through
$\Gamma_K$ by using the word problem algorithm to check
that each $c_i$ is sent to the identity element of \pim.
If the map does not factor, then $\rho$ does not have an inverse and therefore
must have nontrivial kernel.

Otherwise, we have constructed a homomorphism 
\[\phi\co\langle T_1,\ldots,T_s|c_1,\ldots,c_r\rangle \longrightarrow \pim \]
$\phi$ from $\langle T_1,\ldots,T_s|c_1,\ldots,c_r\rangle$ to \pim and 
we have
$\rho\circ\phi=\mathrm{Id}_{\Gamma_K}$ by construction.  To
determine whether $\phi\circ\rho=\mathrm{Id}_{\pim}$ the procedure checks
that for each generator $\phi\circ\rho(g_i)=g_i$, again using the word
problem algorithm.  If any $\phi\circ\rho(g_i)\neq g_i$, then
$\rho((\phi\circ\rho)(g_i)g_i^{-1})=\rho(g_i)\rho(g_i^{-1})=I$ implies that
$\rho$ has kernel.  Otherwise,
$\rho$ is an isomorphism onto $\Gamma_K$.
\end{proof}
\begin{theorem}\label{th:yes}
There exists an algorithm which takes as input a 
triangulated orientable closed 
three-manifold $M$, a
solution to the word problem in $\pim$ and a representation $\rho\co\pim\to
SL_2\Qbar$, and decides if the representation is discrete and
faithful with torsion--free image.
\end{theorem}
\begin{proof}
Let $G=\{g_1,\ldots,g_n\}$ from the
presentation (\ref{eqn:pim}) in Section \ref{subsec:reps}.  We begin
systematically listing elements of the free group $\mathcal{F}(G)$ in such 
a way that an element and its inverse are always added at the same time.  At
each stage of the listing, we have some finite list $W$ of words and a list
$P=\rho(W)\{\pm I\}$ of matrices in $SL_2\Qbar$.  
From this data we attempt to build a fundamental domain.  

Taken together, Lemmas \ref{lemma:spherecoords}, \ref{lemma:NE}, \ref{lemma:K},
\ref{lemma:twosphere}, \ref{lemma:identifiable}, \ref{lemma:fundamental} and
\ref{lemma:gammakisgamma}
give an algorithm which either constructs or fails to construct a 
fundamental domain for the action of $\Gamma$, given
some finite $P\subset\Gamma$ which is closed under inverses.  This algorithm is
outlined in Figure \ref{figure:construct}.
\begin{figure}[ht!]
\begin{center}
\begin{picture}(0,0)%
\includegraphics{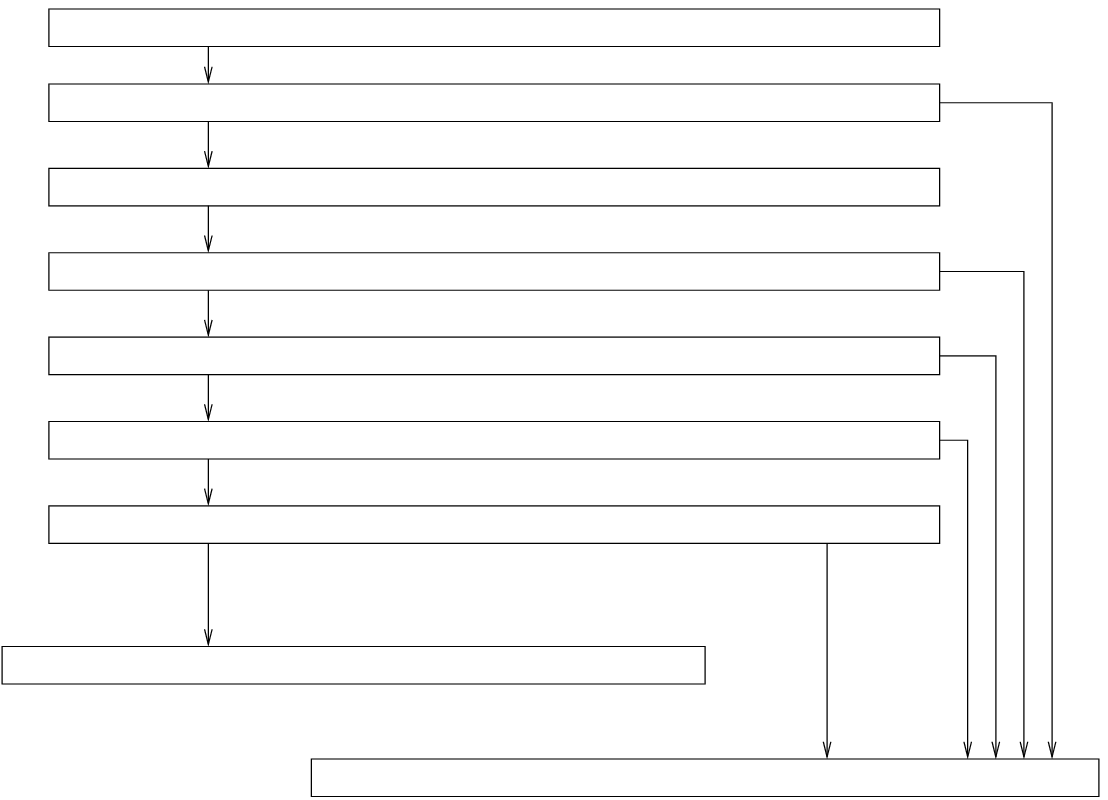}%
\end{picture}%
\setlength{\unitlength}{2368sp}%
\begingroup\makeatletter\ifx\SetFigFont\undefined%
\gdef\SetFigFont#1#2#3#4#5{%
  \reset@font\fontsize{#1}{#2pt}%
  \fontfamily{#3}\fontseries{#4}\fontshape{#5}%
  \selectfont}%
\fi\endgroup%
\begin{picture}(8799,6324)(2089,-7423)
\put(2551,-2611){\makebox(0,0)[lb]{\smash{\SetFigFont{7}{8.4}{\familydefault}{\mddefault}{\updefault}
Construct $K$. (Lemma \ref{lemma:K})
}}}
\put(2551,-3286){\makebox(0,0)[lb]{\smash{\SetFigFont{7}{8.4}{\familydefault}{\mddefault}{\updefault}
Is $K\cong S^2$? (Lemma \ref{lemma:twosphere})
}}}
\put(2551,-3961){\makebox(0,0)[lb]{\smash{\SetFigFont{7}{8.4}{\familydefault}{\mddefault}{\updefault}
Can the faces of $K$ be identified by isometries in $P$? (Lemma \ref{lemma:identifiable})
}}}
\put(2551,-4636){\makebox(0,0)[lb]{\smash{\SetFigFont{7}{8.4}{\familydefault}{\mddefault}{\updefault}
Are the hypotheses of the Poincar\'e Polyhedron Satisfied? (Lemma \ref{lemma:fundamental})
}}}
\put(2551,-5311){\makebox(0,0)[lb]{\smash{\SetFigFont{7}{8.4}{\familydefault}{\mddefault}{\updefault}
Do the face pairings generate $\Gamma$? (Lemma \ref{lemma:gammakisgamma})
}}}
\put(2551,-1936){\makebox(0,0)[lb]{\smash{\SetFigFont{7}{8.4}{\familydefault}{\mddefault}{\updefault}
Does $P$ satisfy condition $NE$?(Lemma \ref{lemma:NE})
}}}
\put(2551,-1336){\makebox(0,0)[lb]{\smash{\SetFigFont{7}{8.4}{\familydefault}{\mddefault}{\updefault}
$P\subset\Gamma$ a finite set closed under inverses
}}}
\put(3901,-2311){\makebox(0,0)[lb]{\smash{\SetFigFont{7}{8.4}{\familydefault}{\mddefault}{\updefault}
YES
}}}
\put(3901,-3586){\makebox(0,0)[lb]{\smash{\SetFigFont{7}{8.4}{\familydefault}{\mddefault}{\updefault}
YES
}}}
\put(3901,-4336){\makebox(0,0)[lb]{\smash{\SetFigFont{7}{8.4}{\familydefault}{\mddefault}{\updefault}
YES
}}}
\put(3901,-4936){\makebox(0,0)[lb]{\smash{\SetFigFont{7}{8.4}{\familydefault}{\mddefault}{\updefault}
YES
}}}
\put(3901,-5761){\makebox(0,0)[lb]{\smash{\SetFigFont{7}{8.4}{\familydefault}{\mddefault}{\updefault}
YES
}}}
\put(9826,-1786){\makebox(0,0)[lb]{\smash{\SetFigFont{7}{8.4}{\familydefault}{\mddefault}{\updefault}
NO
}}}
\put(9826,-3136){\makebox(0,0)[lb]{\smash{\SetFigFont{7}{8.4}{\familydefault}{\mddefault}{\updefault}
NO
}}}
\put(9826,-3811){\makebox(0,0)[lb]{\smash{\SetFigFont{7}{8.4}{\familydefault}{\mddefault}{\updefault}
NO
}}}
\put(9751,-4411){\makebox(0,0)[lb]{\smash{\SetFigFont{7}{8.4}{\familydefault}{\mddefault}{\updefault}
NO
}}}
\put(8776,-5761){\makebox(0,0)[lb]{\smash{\SetFigFont{7}{8.4}{\familydefault}{\mddefault}{\updefault}
NO
}}}
\put(4654,-7336){\makebox(0,0)[lb]{\smash{\SetFigFont{8}{9.6}{\familydefault}{\mddefault}{\updefault}
We cannot construct a fundamental domain for $\Gamma$ from $P$.
}}}
\put(2178,-6436){\makebox(0,0)[lb]{\smash{\SetFigFont{8}{9.6}{\familydefault}{\mddefault}{\updefault}
We have constructed a fundamental domain for $\Gamma$.
}}}
\end{picture}
\caption{An algorithm which tries to construct a fundamental domain for
$\Gamma<SL_2\Qbar$ from a
finite subset of $\Gamma$.  
$K$ is the ``compact part'' of the boundary of $D$, if
$D$ is the region outside all the isometric spheres of elements of $P$.}
\label{figure:construct}
\end{center}
\end{figure}

If this algorithm finds a fundamental domain for the action of $\Gamma$, then it
follows that $\Gamma$ acts as a discrete fixed point free group of isometries
of hyperbolic space.  According to Lemma \ref{lemma:faithful}, there is then an
algorithm to
determine whether $\rho$ is faithful.

Let $\rho\co\pim\to SL_2\Qbar$ be a discrete faithful representation
with torsion--free image $\Gamma$.  We must show that the algorithm we have
described eventually discovers a fundamental domain for the action of $\Gamma$.
Note first that any subset of $\Gamma$ satisfies condition $NE$.  For suppose
that $A$ and $B$ are hyperbolics in $\Gamma$ which
have the same isometric sphere $S_A=S_B$.  Then $AB^{-1}$ fixes \zo, and
so must be trivial or elliptic.  But $\Gamma$ contains no elliptics, so we must
have $A=B$.

The intersection of the exteriors of \emph{all} the isometric spheres of
elements of $\Gamma$ is the Dirichlet domain for $\Gamma$ centered at \zo.  
Theorem \ref{th:homhom} implies that the quotient
of hyperbolic space by $\Gamma$ is homeomorphic to $M$, so it is a closed
3--manifold.  In particular, $\Gamma$ is geometrically finite, so the Dirichlet
domain centered at \zo\ has finitely many sides.  
Each of these sides is part of the isometric sphere of a unique loxodromic
element of $\Gamma$.  There is therefore a finite subset of $\Gamma$ containing
all the matrices whose isometric spheres form faces of the Dirichlet domain.
Eventually our list $P$ will contain all these matrices, and we suppose for the
remainder of the proof that it does.  But now the complex $K$ constructed in
Lemma \ref{lemma:K} is precisely the boundary of the Dirichlet domain centered
at \zo, which is homeomorphic to a sphere.  Face pairings are given by the
isometries whose isometric spheres form the faces, and the algorithm of Lemma
\ref{lemma:identifiable} will verify this.  The Dirichlet domain
satisfies the hypotheses of the Poincar\'e Polyhedron Theorem, as will be
verified by the algorithm of \ref{lemma:fundamental}.  The face identifications
of the Dirichlet domain for an action generate that action, so the algorithm of
\ref{lemma:gammakisgamma} will correctly determine that $\Gamma_K = \pim$.  Thus
if $P$ is large enough, the algorithm of Figure~\ref{figure:construct} will
eventually produce a fundamental domain for the action of $\Gamma$.
\end{proof}

\section{Main theorems}\label{sec:structure}

The main theorem of the paper now follows easily from Theorem \ref{th:no} and
\ref{th:yes}.
\begin{theorem}\label{th:main}
There exists an algorithm which will, given a triangulated orientable closed 
three-manifold $M$ and a
solution to the word problem in $\pim$, decide whether or not $M$ has a
hyperbolic structure.
\end{theorem}
\begin{proof}
By Theorem \ref{th:homhom} $M$ admits a hyperbolic structure if and only if
\begin{enumerate}
\item $M$ is irreducible, and
\item there is an injective homomorphism $\rho\co\pim \to SL_2\C$ 
so that $\rho(\pim)$
acts freely and cocompactly on $\H^3$.
\end{enumerate}
By Lemma \ref{lemma:irreducible} we may suppose that $M$ is irreducible.  
Lemma
\ref{lemma:reps} gives a way to construct a finite list $\mathcal{L}$ of 
representations of
\pim\ into $SL_2\Qbar \subset SL_2\C$, with the property that if a representation
satisfying condition 2 exists, then it is in $\mathcal{L}$. 

Theorems \ref{th:yes} and \ref{th:no} together give a way to decide whether a
particular representation from $\mathcal{L}$ satisfies condition 2 or not.  The
manifold $M$ admits a hyperbolic structure if and only if
some representation in $\mathcal{L}$ satisfies the condition. 
\end{proof}

\begin{theorem}\label{th:structure}
If $M$ is a triangulated closed orientable $3$--manifold which has a
hyperbolic structure, then the hyperbolic structure is algorithmically
constructible.
\end{theorem}
\begin{proof}
Note that by Remark \ref{remark:mapfollows}, it is sufficient to find a
triangulated hyperbolic manifold which is homeomorphic to $M$.
If it is known that $M$ is hyperbolic, it follows that the fundamental group is
automatic.  There are algorithms to find the automatic structure once one
is known to exist \cite{ep:word}.  This automatic structure immediately provides
a procedure for solving the word problem.  We can then use the algorithms of
Sections \ref{sec:alg} and \ref{sec:yes} to find a particular representation of
\pim\ which is discrete and faithful, and a 
fundamental domain for the action of \pim\ on hyperbolic space
through that representation.  The fundamental polyhedron with
face--identifications give a cell--structure for $N = \H^3 / \pim$ which is
easily turned into a triangulation.  
\end{proof}

\end{document}